\documentclass{nyj}

\RequirePackage{amsmath}
\RequirePackage{amsopn}
\RequirePackage{amsfonts}
\RequirePackage{amsthm}

\usepackage[all]{xy}

\title{Tensor products with bounded continuous functions}

\author{Dana P. Williams} 
\address{Department of Mathematics \\ 
  Bradley Hall \\ Hanover, NH 03755-3551 \\ USA}
\email{dana.williams@dartmouth.edu, 
http://www.math.dartmouth.edu/\string~dana}



\CompileMatrices
\allowdisplaybreaks[2]

\newtheorem{thm}{Theorem}

\newtheorem{cor}[thm]{Corollary}
\newtheorem{lemma}[thm]{Lemma}
\newtheorem{prop}[thm]{Proposition}

\theoremstyle{definition}

\theoremstyle{remark}
\newtheorem{remark}[thm]{Remark}

\def\mathcs{C^{*}}
\newcommand{\cs}{\ensuremath{\mathcs}}

\DeclareMathSymbol{\rtimes}{\mathbin}{AMSb}{"6F}

\def\R{\mathbf{R}}
\def\C{\mathbf{C}}

\DeclareMathOperator*{\supp}{supp}

\def\nbhd{neighborhood}

\def\set#1{\{\,#1\,\}}
\let\tensor=\otimes

\makeatletter
\def\labelenumi{\textnormal{(\@alph\c@enumi)}}
\def\theenumi{\@alph \c@enumi}
\newcount\charno
\def\alphapart#1{\charno=96
\advance\charno by#1\char\charno}

\makeatother
\newcommand{\cb}{C^{b}}
\newcommand{\atensor}{\odot}
\newcommand{\cbx}{\cb(X)}
\newcommand{\cbxa}{\cb(X,A)}
\newcommand{\iotaone}{\iota_{1}}
\newcommand{\iotatwo}{\iota_{2}}
\newcommand{\U}{\mathcal{U}}
\newcommand{\V}{\mathcal{V}}
\newcommand{\bx}{\beta X}
\newcommand{\cby}{\cb(Y)}
\newcommand{\cbxcy}{\cb(X\times Y)}
\renewcommand{\phi}{\varphi}

\begin{document}
\begin{abstract}
  We study that natural inclusions $\cbx\tensor A$ into $\cbxa$ and
  $\cb\bigl(X,\cby\bigr)$ into $\cbxcy$.  In particular, excepting
  trivial cases, both these maps are isomorphisms only when $X$ and
  $Y$ are pseudocompact.  This implies a result of Glicksberg showing
  that the Stone-\v Cech compactificiation $\beta(X\times Y)$ is
  naturally identified with $\bx\times \beta Y$ if and only if $X$ and
  $Y$ are pseudocompact.
\end{abstract}
\maketitle
\tableofcontents

Suppose that $X$ is a locally compact Hausdorff space and that $A$ is
a \cs-algebra.  The collection $\cb(X,A)$ of bounded continuous
$A$-valued functions on $X$ is a \cs-algebra with respect to the
supremum norm.  (When $A=\C$, we write simply $\cb(X)$).  Elements in
the algebraic tensor product $\cb(X)\atensor A$ will always be viewed
as functions in $\cb(X,A)$, and the supremum norm on $\cb(X,A)$
restricts to a \cs-norm on $\cbx\atensor A$.  Thus we obtain an
injection $\iotaone$ of the completion $\cbx\tensor A$ into $\cbxa$:
\begin{equation*}
  \iotaone:\cbx\tensor A\hookrightarrow \cbxa,
\end{equation*}
and we can identify $\cbx\tensor A$ with a subalgebra of $\cbxa$.  It
is one of the fundamental examples in the theory that $\iotaone$ is an
isomorphism in the case that $X$ is compact
\cite[Proposition~B.16]{rw:morita}, and we want to investigate the
general case here.  Our main result identifies the range of $\iotaone$
as those functions in $\cbxa$ whose range has compact
closure.  As a consequence we show --- provided $A$ is infinite
dimensional --- that $\iotaone$ is an isomorphism if and only if $X$
has the property that every continuous function on $X$ is bounded.
Such spaces are called \emph{pseudocompact}.  Since paracompact
pseudocompact spaces are compact, it was tempting to make a blanket
assumption of paracompactness.  However, the arguments here use
properties naturally associated to pseudocompactness rather than
compactness, so it seemed worth the little bit of extra effort to
include the general results.  (Nevertheless, I have tried to organize
the paper so that the niceties of noncompact pseudocompact spaces come
at the end.)  In fact, pseudocompactness arises again when we consider
the case where $A=\cby$ for a locally compact Hausdorff space $Y$.
Then we are led to address the properties of the natural inclusion
\begin{equation*}
  \iotatwo:\cb\bigl(X,\cby\bigr) \hookrightarrow \cbxcy.
\end{equation*}
In general, this map is an isomorphism when $\cby$ is given the strict
topology viewed as the multiplier algebra of $C_{0}(Y)$
\cite[Corollary~3.4]{apt:jfa73}.  Here, however, we are interested in
the norm topology, and in Theorem~\ref{thm-main} we show that
$\iotatwo$ is an isomorphism if and --- assuming that $X$ is both
infinite and pseudocompact --- only if $Y$ is pseudocompact.
Combining these observations about $\iotaone$ and $\iotatwo$ yields a
well known result about products of Stone-\v Cech compactifications
originally due to Glicksberg \cite[Theorem~1]{gli:tams59} with
simplified proofs given by Frol\'\i k \cite{fro:cmj60} and Todd
\cite{tod:cmb71}.  Recall that the if $X$ is a locally compact, then,
following \cite[\S XI.8.2]{dug:topology} for example, the Stone-\v
Cech compactification of $X$ is a compact Hausdorff space $\bx$
together with a homeomorphism $i^{X}$ of $X$ onto a dense open subset
of $\bx$ so that $(\bx,i^{X})$ has the \emph{extension property}:
given any continuous map $f$ of $X$ into a compact Hausdorff space
$Y$, there is a continuous map $f^{*}:\bx\to Y$ such that
$f=f^{*}\circ i^{X}$.  (The pair $(\bx,i^{X})$ is unique up to the
natural notion of equivalence.)  
In particular, if $Y$ is locally compact Hausdorff, the
extension property gives a surjection $\phi:\beta(X\times
Y)\to\bx\times \beta Y$.  Our results can be used to show that $\phi$
is a homeomorphism if and only if $X$ and $Y$ are pseudocompact (which
is the special case of Glicksberg's result alluded to above).

\section{Main Results}
\label{sec:main-theorem}

Recall that a subset of a topological space is called
\emph{precompact} if its closure is compact.
\begin{thm}
  \label{thm-one}
  If $X$ is a locally compact Hausdorff space and if $A$ is a
  \cs-algebra, then $f\in\cbxa$ is in $\cbx\tensor A$ if and only if
  the range of $f$, $R(f):=\set{f(x):x\in X}$,
  is precompact.
\end{thm}
\begin{proof}
  Suppose that $f\in\cbx\tensor A$.  To show that $R(f)$ is
  precompact, it will suffice to see that $R(f)$ is totally bounded.
  So we fix $\epsilon>0$ and try to show that $R(f)$ can be covered by
  finitely any $\epsilon$-balls.  Choose $\sum_{i=1}^{n} f_{i}\tensor
  a_{i} \in \cbx\atensor A$ so that
  \begin{equation*}
    \Bigl\|f- \sum_{i=1}^{n} f_{i}\tensor
    a_{i}\Bigr\|_{\infty}<\frac\epsilon2. 
  \end{equation*}
  Since
\begin{equation*}
  \Bigl\|\Bigl(\sum_{i=1}^{n} f_{i}\tensor a_{i}\Bigr)(x) -
  \Bigl(\sum_{i=1}^{n} f_{i}\tensor a_{i}\Bigr)(y)\Bigr\| \le
  \sum_{i=1}^{n} \bigl|f_{i}(x)-f_{i}(y)\bigr|\|a_{i}\|,
\end{equation*}
and since the image of each $f\in\cbx$ is bounded, we can find points
$x_{1},\dots,x_{r}\in X$ such that the open sets
\begin{equation*}
  U_{k}:= \Bigl\{\,x\in X: \Bigl\|\Bigl(\sum_{i=1}^{n}f_{i}\tensor
  a_{i}\Bigr)(x) - \Bigl(\sum_{i=1}^{n}f_{i}\tensor a_{i}\Bigr)(x_{k})
  \Big\|<\frac\epsilon2\,\Bigr\}
\end{equation*}
cover $X$.  For convenience, let $b_{k}:=
\bigl(\sum_{i=1}^{n}f_{i}\tensor a_{i}\bigr)(x_{k})$, and let
\begin{equation*}
  B_{\epsilon}(b_{k}) =\set{a\in A: \|a-b_{k}\|<\epsilon},
\end{equation*}
be the ball of radius $\epsilon$ at $b_{k}$.  Now if $x\in U_{k}$,
then
\begin{align*}
  \Bigl\| f(x)-\Bigl(\sum_{i=1}^{n} f_{i}\tensor
  a_{i}\Bigr)(x_{k})\Bigr\| &\le \Bigl\| f(x) - \Bigl(\sum_{i=1}^{n}
  f_{i} \tensor a_{i}\Bigr)(x)\Bigr\| + \Bigl\| \Bigl(\sum_{i=1}^{n}
  f_{i} \tensor a_{i}\Bigr)(x) -b_{k}\Bigr\| \\
  &<\frac\epsilon2+\frac\epsilon2=\epsilon.
\end{align*}
Since the $U_{k}$ cover $X$, it follows that
\begin{equation*}
  R(f)\subset\bigcup_{k=1}^{r} B_{\epsilon}(b_{k}).
\end{equation*}
Thus $R(f)$ is totally bounded, and $\overline{R(f)}$ must be compact.

Now assume that $\overline{R(f)}$ is compact and therefore totally
bounded.  Thus given $\epsilon>0$, there are elements
$\set{b_{k}}_{k=1}^{r} \subset A$ such that
\begin{equation*}
  \overline{R(f)} \subset \bigcup_{k=1}^{r} B_{\epsilon}(b_{k}).
\end{equation*}
Let $U_{k}:=f^{-1}\bigl(B_{\epsilon}(b_{k})\bigr)$.  Then
$\U=\set{U_{k}}_{k=1}^{r}$ is an open cover of $X$.  First, we suppose
that there is a partition of unity on $X$ subordinate to $\U$.  That
is, we assume that there are functions $f_{k}\in\cbx$ with $0\le
f_{k}\le 1$, $\supp f_{k}\subset U_{k}$ and $\sum_{k=1}^{r}f_{k}(x)=1$
for all $x\in X$.  Then $\sum_{k=1}^{r} f_{k}\tensor
b_{k}\in\cbx\atensor A$, and for all $x\in X$, we have
\begin{align*}
  \Bigl\|\Bigl(\sum_{k=1}^{r} f_{k}\tensor b_{k}\Bigr)(x)-f(x)\Bigr\|
  &= \Bigl\| \sum_{k=1}^{r} f_{k}(x)b_{k} -\sum_{k=1}^{r} f_{k}(x)f(x)
  \Bigr\| \\
  &\le \sum_{k=1}^{n} f_{k}(x)\|b_{k}-f(x)\| \\
  &\le \epsilon\sum_{k=1}^{n} f_{k}(x) 
  = \epsilon.
\end{align*}
Thus, if such a partition of unity exists, then $f$ belongs to the
closure of $\cbx\atensor A$ and $f\in\cbx\tensor A$ as required. If
$X$ were paracompact, then given any finite cover $\U$, there is a
partition of unity $\set{f_{k}}_{k=1}^{r}$ on $X$ subordinate to $\U$
\cite[Proposition~4.34]{rw:morita}. Since, out of stubbornness, $X$
is not assumed to be paracompact, we will have to take advantage
of the special nature of the covers $\U$ involved and the extension
property of the Stone-\v Cech compactification $\bx$ of $X$.  Since
$\overline{R(f)}$ is compact and $f:X\to \overline{R(f)}$ is
continuous, there is a continuous function $F:\bx\to \overline{R(f)}$
such that $F\circ i^{X}=f$.  Let
$V_{k}:=F^{-1}\bigl(B_{\epsilon}(b_{k})\bigr)$.  Then $\V=\set{V_{k}}$ is
a finite open cover of $\bx$, and there is a partition of unity
$\set{\phi_{k}}$ on $\bx$ subordinate to $\V$.  Since $U_{k} =
(i^{X})^{-1}(V_{k})$, it follows that if $f_{k}:= \phi_{k}\circ
i^{X}$, then $\set{f_{k}}$ is a partition of unity on $X$ subordinate
to $\U$.  This completes the proof.
\end{proof}

\begin{cor}
  \label{cor-main}
  Suppose that $X$ is a locally compact Hausdorff space and that $A$
  is a \cs-algebra.  If $X$ is pseudocompact, then $R(f)$ is compact
  for each $f\in\cbxa$, and
  \begin{equation}\label{eq:2}
    \cbx\tensor A=\cbxa
  \end{equation}
  (that is, $\iotaone$ is an isomorphism).  Conversely, if $A$ is not
  finite dimensional, then \eqref{eq:2} holds only if $X$ is
  pseudocompact.  If $A$ is finite dimensional, then \eqref{eq:2}
  always holds.
\end{cor}
\begin{proof}
  Suppose that $X$ is pseudocompact.
  It is not hard to see that the continuous image of a paracompact
  space in a metric space is compact (Corollary~\ref{cor-ctspc}). Thus
  $R(f)$ is compact for each $f\in \cbxa$ and \eqref{eq:2} follows
  from Theorem~\ref{thm-one}.
  
  Now assume that $A$ is infinite dimensional and that $X$ is not
  pseudocompact.  It follows from
  \cite[Theorem~1.23]{rud:functional}, that the unit ball of $A$ is
  not compact.  Thus there is a sequence $\set{a_{k}}_{k=1}^{\infty}$
  of elements of $A$ of norm at most one such that
  $\set{a_{k}}_{k=1}^{\infty}$ is not totally bounded.  Since $X$ is
  not pseudocompact, there are precompact open sets $U_{n}$ in $X$
  whose closures are locally finite and pairwise
  disjoint (Lemma~\ref{lem-notpc}). 
  Fix $x_{n}\in U_{n}$.  Then by
  Urysohn's Lemma (cf. \cite[Proposition~1.7.5]{ped:analysis}), there
  are 
  $f_{n}\in C_{c}(X)$ such that $0\le f_{n}\le 1$, $f_{n}(x_{n})=1$
  and $\supp f_{n}\subset U_{n}$.  Since $\set{\overline{U_{n}}}$ is
  locally finite,
\begin{equation*}
  f(x):=\sum_{n=1}^{\infty} f_{n}(x)a_{n}
\end{equation*}
defines a function $f\in\cbxa$ whose range contains
$\set{a_{n}}_{n=1}^{\infty}$.  Thus $f$ is not in the range of
$\iotaone$, and \eqref{eq:2} does not hold.

Since every finite dimensional \cs-algebra is a direct sum of matrix
algebras \cite[Theorem~6.3.8]{mur:cs-algebras} and since it is easy
to see that $\cb(X,M_{n})\cong M_{n}\bigl(\cbx\bigr)$, the last
assertion is an easy consequence another standard example in the
theory of tensor products: $M_{n}\bigl(\cbx\bigr) \cong \cbx\tensor
M_{n} $ \cite[Proposition~B.18]{rw:morita}.
\end{proof}

Now we want consider the case where $A$ is of the form
$\cb(Y)$ for some locally compact Hausdorff space $Y$, and investigate
the natural inclusion~$\iotatwo$.
\begin{thm}
  \label{thm-main}
  Suppose that $X$ and $Y$ are locally compact Hausdorff spaces.  If
  $Y$ is pseudocompact, then the natural inclusion
  $\iotatwo:\cb\bigl(X,\cby\bigr) \hookrightarrow \cbxcy$ is an
  isomorphism.  If $X$ is pseudocompact and infinite, then $\iotatwo$
  is an isomorphism only if $Y$ is pseudocompact.
\end{thm}

We'll postpone the proof of until \S\ref{sec:proof-proposition} so
that we can see how Glicksberg's result follows from
Theorems~\ref{thm-one}~and \ref{thm-main}.

\section{Glicksberg's Theorem}
\label{sec:glicksbergs-theorem}

  If $Z$ is
  compact Hausdorff and if $h:X\to Z$ is a homeomorphism of $X$ onto a
  dense subset of $Z$, then $h(X)$ is open in $Z$
  \cite[\S XI.8.3]{dug:topology}, and the extension property of $\bx$
  implies that there is a unique continuous
  surjection $\phi:\bx\to Z$ such that $\phi\circ i^{X}=h$.
Equivalently, we get a commutative diagram of algebra
homomorphisms 
\begin{equation}\label{eq:10}
\begin{split}
  \xymatrix{&C(\bx)\\ C_{0}(X) \ar[ru]^{i^{X}_{*}} \ar[r]_{h_{*}} &
  C(Z) \ar@{.>}[u]_{\phi^{*}},}
\end{split}
\end{equation}
where, for example,
\begin{equation*}
  h_{*}(f)(z)=
  \begin{cases}
    0&\text{if $y\notin h(X)$, and}\\
    f(x)&\text{if $z=f(x)$ for some $x\in X$,}
  \end{cases}
\end{equation*}
and $\phi^{*}(f)(z)=f\bigl(\phi(z)\bigr)$.  Note that $\phi^{*}$ is
the unique homomorphism making \eqref{eq:10} commute.  In particular,
if $X$ and $Y$ are locally compact spaces, then we have a commutative
diagram
\begin{equation*}
  \xymatrix{&\beta(X\times Y) \ar[d]^{\phi} \\ X\times Y
  \ar[ur]^{i^{X\times Y}} \ar[r]_-{i^{X}\times i^{Y}}& \bx\times\beta Y,}
\end{equation*}
and it is natural to ask when $\phi$ is a homeomorphism so that we can
identify $(\bx\times\beta Y,i^{X}\times i^{Y})$ with
$\bigl(\beta(X\times Y),i^{X\times Y}\bigr)$.  Since $\phi$, and hence
$\phi^{*}$, is unique and since the extension property of $\bx$ easily
implies $\cbx\cong C(\bx)$, we can find $\phi$ by combining
the following natural maps:
\begin{equation*}
  \xymatrix{
&C\bigl(\beta(X\times Y)\bigr) \ar[rr]^{\cong} &&\cb(X\times Y) \\
C_{0}(X\times Y) \ar[dr]_{(i^{X}\times i^{Y})_{*}} \ar[ur]^{i^{X\times
  Y}_{*}} &&&\cb\bigl(X,\cby\bigr) \ar[u]_{\iotatwo} \\
&C(\bx\times \beta Y)\ar@{.>}[uu]_{\phi^{*}}
  \ar[r]_-{\cong} & C(\bx)\tensor 
  C(\beta Y)\ar[r]_{\cong} &\cbx\tensor \cby \ar[u]_{\iotaone}.}
\end{equation*}
Thus $\phi^{*}$ is an isomorphism exactly when both $\iotaone$ and
$\iotatwo$ are surjective.  Thus if both $X$ and $Y$ are infinite, so
that, for example, $A=\cby$ is infinite dimensional, then $\iotaone$ is
an isomorphism if and only if $X$ is pseudocompact
(Corollary~\ref{cor-main}), and if $X$ is pseudocompact, then 
$\iotatwo$ is an isomorphism if and
only if $Y$ is pseudocompact (Theorem~\ref{thm-main}).  Since
$\phi^{*}$ is an isomorphism exactly when $\phi$ is a homeomorphism,
we have proved a special case of Glicksberg's result
\cite[Theorem~1]{gli:tams59}.

\begin{thm}[Glicksberg]
  \label{thm-glicksberg}
  Suppose that $X$ and $Y$ are infinite locally compact spaces.  Then
  the natural map $\phi:\beta(X\times Y) \to \bx\times\beta Y$ is a
  homeomorphism if and only if both $X$ and $Y$ are pseudocompact.
\end{thm}

Of course, Glicksberg considered arbitrary products and only assumed
that that $X$ and $Y$ are completely regular.  We have dispensed with
completely regular spaces out of prejudice, and the extension to
arbitrary products is not difficult and is discussed in \cite{tod:cmb71}.

\begin{remark}
  \label{rem-path}
  Note that even if one of $X$ and $Y$ fails to be pseudocompact,
  Theorem~\ref{thm-glicksberg} does not preclude the possibility that
  $\beta(X\times Y)$ and $\bx\times \beta Y$ are homeomorphic.  It
  only asserts that the natural, and only useful, way to identify them
  fails.  See \cite[\S6]{gli:tams59} for further thoughts on this.
\end{remark}

\section{Pseudocompact spaces}
\label{sec:pseudocompact-spaces}

Before turning to the proof of Theorem~\ref{thm-main}, some
preliminaries on pseudocompact spaces are in order.  Of course,
everything sketched here is well known, but distilling the specifics
from the literature could be tedious.
\begin{lemma}
  \label{lem-notpc}
  If $X$ is locally compact and \emph{not} pseudocompact, then there is a
  sequence $\set{V_{n}}$ of nonempty open sets in $X$ such that each
  point in $X$ meets at most one $\overline{V_{n}}$.
\end{lemma}
\begin{proof}
  Since $X$ is not pseudocompact, there is a unbounded
  continuous nonnegative real-valued function $f$ on $X$.  Thus there is a
  sequence $\set{x_{n}}\subset X$ such that $f(x_{n+1})> f(x_{n})+1$ for
  all $n$.
  By composing with a continuous piecewise linear function on $\R$, we
  may as well assume that $f(x_{n})=n$.  Then we can let
  $V_{n}:=f^{-1}\bigl((n-\frac14,n+\frac14)\bigr)$.
\end{proof}

Recall that a family $\set{U_{i}}$ of subsets of $X$ is called
\emph{locally finite} if every point in $X$ has a \nbhd{} meeting at
most finitely many $U_{i}$.

\begin{prop}
  \label{prop-pc}
  A locally compact space $X$ is pseudocompact if and only if every
  countable  
  locally finite collection of nonempty open sets in $X$ is finite.
\end{prop}
\begin{proof}
  Let $\set{U_{n}}$ be locally finite sequence of nonempty open sets
  in $X$.  Fix $x_{n}\in U_{n}$.  By Urysohn's Lemma, there is 
a continuous
  function $f_{n}$ on $X$ such that $0\le f_{n}\le 1$,
  $f_{n}(x_{n})=n$ and $\supp f_{n } \subset U_{n}$.  Since
  $\set{U_{n}}$ is locally finite, $f=\sum f_{n}$ is continuous (and
  unbounded) on $X$.  Thus $X$ is not pseudocompact.  This, together
  with Lemma~\ref{lem-notpc}, establishes the result.
\end{proof}

Since, by definition, every open cover of a paracompact space has a
locally finite subcover, it follows immediately from
Proposition~\ref{prop-pc} that a pseudocompact paracompact space is
compact.  Since metric spaces are paracompact and since the inverse
image of a locally finite cover is locally finite, we
obtain the following easy corollary of Proposition~\ref{prop-pc}.
\begin{cor}
  \label{cor-ctspc}
  If $X$ is pseudocompact and if $f:X\to Y$ is a continuous function
  from $X$ into a metric space $Y$, then $f(X)$ is compact.
\end{cor}

Nevertheless, there certainly are locally compact noncompact
pseudocompact spaces.  For example, if $\omega_{1}$ is the first
uncountable ordinal, then $[0,\omega_{1})$ is a countably compact
locally compact space in the order topology, and countably compact
spaces are easily seen to be pseudocompact.  Moreover, if $\omega$ is
the first countable ordinal, then $[0,\omega_{1}]\times
[0,\omega]\setminus\set{(\omega_{1},\omega)}$ is pseudocompact and not
even countably compact \cite[\S XI.3 Ex.~2]{dug:topology}.  More
generally, $X$ is pseudocompact if and only if the corona set
$\bx\setminus X$ contains no nonempty closed $G_{\delta}$ sets
\cite[p.~370]{gli:tams59}.  This observation --- together with an old
result of \v Cech showing that closed $G_{\delta}$ sets in the corona
of a noncompact locally compact space always have cardinality at least
that of the continuum \cite[p.~835]{cec:am37} --- has a curious
consequence.
If
$x\in\bx\setminus X$, then $\bx\setminus\set x$ is always
pseudocompact \cite[pp.~380--1]{gli:tams59}.

As it happens, the product of two locally compact pseudocompact spaces
is pseudocompact \cite[Theorem~4(a)]{gli:tams59}.  We give a short
proof of the special case of this we need in
\S\ref{sec:proof-proposition}.  Some care is called for as, in
general, the product of (not necessarily locally compact)
pseudocompact spaces need not be pseudocompact
\cite[p.~245]{dug:topology}.
\begin{lemma}
  \label{lem-added}
  Suppose that $X$ and $Y$ are locally compact Hausdorff spaces with
  $X$ compact and $Y$ pseudocompact.  Then $X\times Y$ is
  pseudocompact.
\end{lemma}
  \begin{proof}
    It will suffice to prove that if $\set{U_{n}}$ and $\set{V_{n}}$
    are infinite sequences of nonempty open sets in $X$ and $Y$,
    respectively, then $\set{U_{n}\times V_{n}}$ is not locally
    finite.  Since $Y$ is pseudocompact, $\set{V_{n}}$ can't be
    locally finite and there is a $y\in Y$ such that every \nbhd{} of
    $y$ meets infinitely many $V_{n}$.  Let
    \begin{equation*}
\Omega:=\set{ (n,W):\text{$n\in \mathbf{N}$, $W$ is a \nbhd{} of
    $y$ and $W\cap V_{n}\not=\emptyset$}}.
\end{equation*}
It is easy to see that $\Omega$ is directed
\cite[Definition~X.1.2]{dug:topology}.  Thus if we choose
$y_{(n,W)}\in W\cap V_{n}$, then $\set{y_{\omega}}_{\omega\in\Omega}$
is a net in $Y$ converging to $y$.  On the other hand, if for each
$(n,W)\in\Omega$ we choose $x_{(n,W)}\in U_{n}$, then, since $X$ is
compact, $\set{x_{\omega}}_{\omega\in\Omega}$ has a subnet
$\set{x_{\omega_{i}}}_{i\in I}$ which converges to some $x\in X$.
Thus if $U$ and $V$ are \nbhd s of $x$ and $y$, respectively, then
$\set{(x_{\omega_{i}},y_{\omega_{i}})}$ is eventually in $U\times V$.
Thus if $\omega_{i}:=(n_{i},W_{i})$, then we eventually have $(U\times
V)\cap (U_{n_{i}}\times V_{n_{i}})\not=\emptyset$.  It follows that
every \nbhd{} of $(x,y)$ eventually meets infinitely many $U_{n}\times
V_{n}$.  Thus $X\times Y$ is pseudocompact.
  \end{proof}

  \begin{remark}
    \label{rem-prod}
    Using the observation that a space $X$ is pseudocompact if and
    only if every bounded continuous function on $X$ attains its
    maximum on $X$, it is not hard to see that
    Theorem~\ref{thm-glicksberg} implies that the product of locally
    compact pseudocompact spaces is pseudocompact
    \cite[Theorem~4(a)]{gli:tams59}.  Let $F\in\cb(X\times Y)$.  Then
    Theorem~\ref{thm-glicksberg} implies that $F$ has a unique
    continuous extension $F^{*}$ to $C(\bx\times \beta Y)$.  Let
    $f^{*}(x):= F(x,\cdot)$.  Since $X$ is pseudocompact, there is a
    $\bar x\in X$ such that $\|f^{*}(\bar x)\|_{\infty}=\sup_{x\in
    X}\|f^{*}(x)\|_{\infty}$.  But since $Y$ is pseudocompact, there
    is a $\bar 
    y\in Y$ such that $\|f(\bar x)\|_{\infty}=f(\bar x)(\bar y)$.  Thus,
    $F$ assumes its maximum at $(\bar x,\bar y)$, and $X\times Y$ is
    pseudocompact. 
  \end{remark}

\section{Proof of Theorem~\ref{thm-main}}
\label{sec:proof-proposition}

If $Y$ is compact, then the surjectivity of $\iotatwo$ is fairly
standard.  However, if $Y$ is only pseudocompact, then the proof of
surjectivity given here depends on a result of Frol{\'{\i}}k's
\cite[Lemma~1.3]{fro:cmj60}.\footnote{This result is also essential
  to both Frol{\'{\i}}k's \cite{fro:cmj60} and Todd's \cite{tod:cmb71}
  proofs of Glicksberg's result.}  Since this
result is the heart of the proof, and since the published version has
some annoying typos, we include the proof here for completeness.

\begin{lemma}[Frol{\'{\i}}k]
  \label{lem-frolik}
  Suppose that $X$ and $Y$ are locally compact Hausdorff spaces with
  $Y$ pseudocompact.  If $f$ is a bounded real-valued function on
  $X\times Y$, then
  \begin{equation*}
    g(x):=\sup_{y\in Y} f(x,y)
  \end{equation*}
  defines a continuous function on $X$.
\end{lemma}
\begin{proof}
  Since continuity is a local property, we may as well assume that $X$
  is compact.  Since $g$ is the supremum of continuous functions on
  $X$, it is easy to see that $g$ is lower semicontinuous; that is,
  $\set{x\in X:g(x)>b}$ is open for all $b\in\R$.  Thus it will
  suffice to see that $g$ is also upper semicontinuous so that
  $\set{x\in X:g(x)<b}$ is open for all $b\in \R$.  If $g$ is not
  upper semicontinuous, then there must be a $x_{0}\in X$ and a
  $\epsilon>0$ such that every \nbhd{} $U$ of $x_{0}$ contains a $x$
  such that $g(x)>g(x_{0})+3\epsilon$.  In particular, there is a
  $(x,y)\in U\times Y$ such that
  \begin{equation}
    \label{eq:5}
    f(x,y)>g(x_{0})+3\epsilon.
  \end{equation}
  We claim that we can choose points $(x_{n},y_{n})$, and \nbhd s
  $V_{n}$ of $y_{n}$, $U_{n}$ of $x_{n}$ and $U_{n}'$ of $x_{0}$ such
  that
\begin{enumerate}
\item $\sup\set{|f(x,y)-f(x',y')|:(x,y),(x',y')\in U_{n}\times
    V_{n}}<\epsilon$, \label{item:1}
\item $\sup\set{|f(x,y)-f(x',y')|:(x,y),(x',y')\in U_{n}'\times
    V_{n}}<\epsilon$, \label{item:2}
\item $U_{n+1}'\subset U_{n}'$, \label{item:4}
\item $U_{n}\subset U_{n-1}'$ provided $n\ge 2$, and\label{item:3}
\item $f(x_{n},y_{n})> g(x_{0})+3\epsilon$.\label{item:5}
\end{enumerate}
For convenience, let $U_{0}'=X$.  Assume we have chosen
$(x_{i},y_{i})$, $V_{i}$, $U_{i}$ and $U_{i}'$ for $i<n$.  Then by
assumption we can choose $(x_{n},y_{n})\in U_{n-1}'\times Y$ such that
(\ref{item:5}) holds.  The existence of $V_{n}$, $U_{n}$ and $U_{n}'$
such that (\ref{item:1})--(\ref{item:3}) hold follows easily from the
continuity of $f$.

Since $X\times Y$ is pseudocompact (Lemma~\ref{lem-added}), there is a
$(\bar x,\bar y)\in X\times Y$ such that every \nbhd{} of $(\bar
x,\bar y)$ meets infinitely many $U_{n}\times Y_{n}$.  Let $U$ and $V$
be \nbhd s of $\bar x$ and $\bar y$, respectively, so that $(x,y)\in
U\times V$ implies
\begin{equation}
  \label{eq:6}
  |f(x,y)-f(\bar x,\bar y)|<\epsilon.
\end{equation}
By assumption, there is a $n$ such that there exists $(x',y')\in
(U_{n}\times V_{n})\cap (U\times V)$.  Using \eqref{eq:6} together with
(\ref{item:1}) and (\ref{item:5}), we have
\begin{equation}
  \label{eq:7}
  f(\bar x,\bar y)> g(x_{0}) + 2\epsilon.
\end{equation}
On the other hand, there is a $m>n$ such that there exists
$(x'',y'')\in (U_{m}\times V_{m})\cap U\times V$.  But
(\ref{item:4})~and (\ref{item:3}) imply that
\begin{equation*}
U_{m}\subset U_{m-1}'\subset U_{n}'.
\end{equation*}
Thus $(x'',y')\in (U_{n}'\times V_{n})\cap (U\times V)$.  Now
\eqref{eq:6} implies $|f(x'',y')-f(\bar x,\bar y)|<\epsilon$, while
(\ref{item:2}) implies $|f(x'',y')-f(x_{0},y')|<\epsilon$.  Thus
\begin{equation*}
  f(\bar x,\bar y)< f(x_{0},y')+\epsilon \le g(x_{0})+\epsilon.
\end{equation*}
This contradicts \eqref{eq:5}, and completes the proof.
\end{proof}

  \begin{proof}[Proof of Theorem~\ref{thm-main}]
    Let $F\in\cbxcy$ and define $f:X\to\cby$ by $f(x)(y):=F(x,y)$.
    Notice that $F$ is in the image of $\iotatwo$ if and only if $f$
    is continuous.  If $Y$ is pseudocompact and $x_{0}\in X$, then
\begin{equation*}
  G(x,y):=|f(x,y)-f(x_{0},y)|
\end{equation*}
is continuous on $X\times Y$.  Frol{\'{\i}}k's Lemma~\ref{lem-frolik}
implies that
\begin{equation*}
  g(x):=\sup_{y\in Y}G(x,y) =\|f(x)-f(x_{0})\|_{\infty}
\end{equation*}
is continuous on $X$.  Since $g(x_{0})=0$, $f$ is continuous at
$x_{0}$.  Since $x_{0}$ was arbitrary, $\iotatwo$ is surjective
whenever $Y$ is pseudocompact.

Now suppose that $X$ is infinite and pseudocompact and that $Y$ is not
pseudocompact.  Let $\set{x_{n}}$ be an infinite sequence in $X$.  A
simple argument (cf. \cite[\S VII.2.4]{dug:topology}) shows that
there are precompact \nbhd s $U_{n}$ of $x_{n}$ such that
$\overline{U_{n}} \cap \overline{U_{m}} =\emptyset $ if $n\not= m$.
Since $Y$ is not pseudocompact, there is a sequence of nonempty
precompact open sets $\set{V_{n}}$ with pairwise disjoint closures 
such that each point in $y$ has a
\nbhd{} meeting at most one $\overline{V_{n}}$ (Lemma~\ref{lem-notpc}).  
It follows that
\begin{equation*}
  \bigcup \overline{U_{n}}\times \overline{V_{n}}
\end{equation*}
is $\sigma$-compact and closed in $X\times Y$. Fix $y_{n}\in V_{n}$.
Then $\set{(x_{n},y_{n})}$ is also closed in $X\times Y$.  Since
$\sigma$-compact locally compact Hausdorff spaces are paracompact,
Urysohn's Lemma implies there is a continuous function $F$ on $\bigcup
\overline{U_{n}}\times \overline{V_{n}}$ such that $0\le F\le 1$,
$F(x_{n},y_{n})=1$ for all $n$ and $\supp F\subset\bigcup U_{n}\times
V_{n}$.  We can extend $F$ to a bounded continuous function on
$X\times Y$ by letting it be identically zero on the complement of
$\bigcup U_{n}\times V_{n}$.

As above, let $f(x)=F(x,\cdot)$.  If $F$ were in the image of
$\iotatwo$, then $x\mapsto f(x)$ would be continuous.  Since $X$ is
pseudocompact, $\set{f(x)}_{x\in X}$ is compact in $\cby$
(Corollary~\ref{cor-ctspc}).  In particular, $\set{f(x_{n})}$ has a
convergent subsequence $\set{f(x_{n_{i}})}$.  Then there is a
$k$ such that $i\ge k$ implies
\begin{equation}
  \label{eq:9}
  \|f(x_{n_{i}})-f(x_{n_{k}})\|_{\infty}<\frac12.
\end{equation}
Of course we can choose $i\ge k$ such that $n_{i}>n_{k}$.  Then
$(x_{n_{i}}, y_{n_{k}})\notin \bigcup U_{n}\times V_{n}$ and
\begin{align*}
  f(x_{n_{i}})(y_{n_{k}}) - f(x_{n_{k}})(y_{n_{k}}) &=
  F(x_{n_{i}},y_{n_{k}}) - F(x_{n_{k}},y_{n_{k}}) \\
  &= 0-1 =-1,
\end{align*}
and this contradicts \eqref{eq:9}.  This completes the proof.
  \end{proof}


\def\noopsort#1{}\def\cprime{$'$}
\ifx\undefined\bysame
\newcommand{\bysame}{\leavevmode\hbox to3em{\hrulefill}\,}
\fi

\end{document}